\documentclass{article}
\usepackage{stmaryrd}
\usepackage{mathrsfs}
\usepackage[centertags]{amsmath}
\usepackage{amsfonts, dsfont}
\usepackage{amssymb}
\usepackage{amsthm}
\usepackage[OT2,T1]{fontenc}
\DeclareSymbolFont{cyrletters}{OT2}{wncyr}{m}{n}
\DeclareMathSymbol{\Sha}{\mathalpha}{cyrletters}{"58}

\input xy
\xyoption{all}

\theoremstyle{plain}
\newtheorem{thm}{Theorem}[section]
\newtheorem{cor}[thm]{Corollary}
\newtheorem{lem}[thm]{Lemma}
\newtheorem{prop}[thm]{Proposition}

\theoremstyle{definition}
\newtheorem{defn}[thm]{Definition}
\newtheorem*{remark}{Remarks}
\newtheorem*{ack}{Acknowledgments}

\newcommand{\bd}{\begin{defn}}
\newcommand{\ed}{\end{defn}}
\newcommand{\bl}{\begin{lem}}
\newcommand{\el}{\end{lem}}
\newcommand{\bp}{\begin{prop}}
\newcommand{\ep}{\end{prop}}
\newcommand{\bt}{\begin{thm}}
\newcommand{\et}{\end{thm}}
\newcommand{\bc}{\begin{cor}}
\newcommand{\ec}{\end{cor}}
\newcommand{\br}{\begin{remark}}
\newcommand{\er}{\end{remark}}
\newcommand{\bdi}{\begin{diagram}}
\newcommand{\edi}{\end{diagram}}
\newcommand{\beq}{\begin{equation}}
\newcommand{\eeq}{\end{equation}}
\newcommand{\ba}{\begin{array}}
\newcommand{\ea}{\end{array}}
\newcommand{\bpf}{\begin{proof}}
\newcommand{\epf}{\end{proof}}

\newcommand{\Q}{\mathds{Q}}
\newcommand{\Zp}{\mathds{Z}_{p}}
\newcommand{\Qp}{\mathds{Q}_{p}}

\newcommand{\al}{\alpha}

\newcommand{\Ga}{\Gamma}

\DeclareMathOperator{\Sel}{Sel} \DeclareMathOperator{\Gal}{Gal}

\newcommand{\cyc}{\mathrm{cyc}}

\newcommand{\ot}{\otimes}
\newcommand{\ilim}{\displaystyle \mathop{\varinjlim}\limits}

\newcommand{\coker}{\mathrm{coker}\,}

\newcommand{\lra}{\longrightarrow}

\newcommand{\ps}[1]{\llbracket #1 \rrbracket}

\begin{document}
\title{On the Euler characteristics of signed Selmer groups}
 \author{ Suman Ahmed\footnote{School of Mathematics and Statistics,
Central China Normal University, Wuhan, 430079, P.R.China.
 E-mail: \texttt{npur.suman@gmail.com}}  \quad
  Meng Fai Lim\footnote{School of Mathematics and Statistics $\&$ Hubei Key Laboratory of Mathematical Sciences,
Central China Normal University, Wuhan, 430079, P.R.China.
 E-mail: \texttt{limmf@mail.ccnu.edu.cn}} }
\date{}
\maketitle

\begin{abstract} \footnotesize
\noindent Let $p$ be an odd prime number, and $E$ an elliptic curve defined over a number field with good reduction at every prime of $F$ above $p$.
In this short note, we compute the Euler characteristics of the signed Selmer groups of $E$ over the cyclotomic $\Zp$-extension. The novelty of our result is that we allow the elliptic curve to have mixed reduction types for primes above $p$ and that we allow mixed signs in the definition of the signed Selmer groups.

\medskip
\noindent Keywords and Phrases: Euler characteristics, signed Selmer groups.

\smallskip
\noindent Mathematics Subject Classification 2010: 11G05, 11R23, 11S25.
\end{abstract}

\section{Introduction}
Let $p$ be an odd prime. Let $F$ be a number field and $E$ an elliptic curve defined over $F$. Suppose for now $E$ has good ordinary reduction at every primes of $F$ above $p$. One can define the $p$-primary Selmer group of $E$ over the cyclotomic $\Zp$-extension $F^{\cyc}$  of $F$. The said Selmer group is conjectured to be cotorsion over $\Zp\ps{\Ga}$ (see \cite{Maz}), where $\Ga=\Gal(F^{\cyc}/F)$. Under this torsionness conjecture, Perrin-Riou \cite{PR} and P. Schneider \cite{Sch} computed the $\Ga$-Euler characteristics of these Selmer groups. The importance of these $\Ga$-Euler characteristics stems from the fact that their values are related to the $p$-part of the algebraic invariants appearing in the formula of the BSD-conjecture which in turn allow one to study the special values of the Hasse-Weil $L$-function of $E$ via the so-called ``Iwasawa Main Conjecture'' (see \cite{CS00, G89, G99, Maz}).

In this paper, we would like to consider the situation where our elliptic curve $E$ may have good supersingular reduction at some primes above $p$. In this case, one usually works with the so-called signed Selmer groups of $E$ in the sense of \cite{KimParity, KimPM,KO,Kob}. Our main result is concerned with computing the Euler characteristics of these signed Selmer groups which we now describe. Suppose that $E$ has good (not necessarily ordinary) reduction at any prime of $F$ lying above $p$. Denote by $S_p^{ord}$ (resp., $S_p^{ss}$) for the set of good ordinary reduction (resp., good supersingular reduction) primes of $E$ above $p$. Suppose further that for each $v\in S_p^{ss}$, one has $F_v=\Qp$ and $a_v = 1 + p - |\tilde{E}_v(\mathbb{F}_p)| = 0$, where $\tilde{E}_v$ denotes the reduction of $E$ at $v$. Our main result is as follows.

\bt \label{main1}
Retain the above settings. Suppose that $\Sel(E/F)$ is finite. Then $\Sel^{\overrightarrow{s}}(E/F^{\cyc})$ is a cotorsion $\Zp\ps{\Ga}$-module and its
 $\Ga$-Euler characteristics is given by
\[ \frac{|\Sha(E/F)(p)|}{|E(F)(p)|^2}\times\prod_v c_v^{(p)}\times\prod_{v\in S_p^{ord}}(d_v^{(p)})^2.\]
 Here $c_v^{(p)}$ is the highest power of $p$ dividing $|E(F_v):E_0(F_v)|$, where $E_0(F_v)$ is the subgroup of $E(F_v)$ consisting of points with nonsingular reduction modulo $v$ and $f_v$ is the residue field of $F_v$, and $d_v^{(p)}$ is the highest power of $p$ dividing $|\tilde{E}_v(f_v)|$.
\et

For the definition of the signed Selmer group $\Sel^{\overrightarrow{s}}(E/F^{\cyc})$, we refer readers to Section \ref{Selmer}. The $\Ga$-Euler characteristics of $\Sel^{\overrightarrow{s}}(E/F^{\cyc})$ is defined to be the quantity
\[ \frac{|H^0(\Ga,\Sel^{\overrightarrow{s}}(E/F^{\cyc}))|}{|H^1(\Ga,\Sel^{\overrightarrow{s}}(E/F^{\cyc}))|}.\]
In the midst of proving Theorem \ref{main1}, we will see that the $\Ga$-Euler characteristics of $\Sel^{\overrightarrow{s}}(E/F^{\cyc})$ makes sense.
When the elliptic curve has good supersingular reduction at all primes above $p$, this formula was first established by Kim in \cite{KimPM}. Our main result improves this prior result in that we allow our elliptic curve to have mixed reduction types for primes above $p$ and that we allow mixed signs in the definition of the signed Selmer groups. The proof of the theorem will be given in Section \ref{Selmer}. In fact, in the said section, we shall consider a slightly more general situation than that stated in this introductory section (see Theorem \ref{main2}). As an application, we establish a result which says that if one of the signed Selmer group vanishes, so do the others (see Corollary \ref{main2cor}).

It would definitely be of interest to be able to provide examples illustrating our theorem. In fact, it is not difficult to obtain examples of elliptic curves with mixed reduction types at primes above $p$ via similar arguments to that in \cite[Proposition 5.4]{G99} or \cite[Lemma 8.19]{Maz}. The problem is that we do not know how to verify the finiteness of $\Sel(E/F)$ in these examples. Until a (nice enough) theory of Euler system has been developed in this mixed reduction context, this does not seem tractable. We do however hope to review this problem in subsequent studies.

Upon the completion of this work, we were informed by Antonio Lei that he and his coauthor have computed the Euler characteristics of the signed Selmer groups over a $\Zp^d$-extension (see \cite{LeiSu}). However in their article, they have worked with elliptic curves with good supersingular reduction at all primes above $p$ and with the same sign in their definition of the signed Selmer groups. There they also required that the prime $p$ to split completely over $F/\Q$. It would be of interest to see if a similar computation can be performed for the situation considered in Section \ref{Selmer} of our paper. One might even contemplate computing these Euler characteristics over a noncommutative $p$-adic extension. We hope to explore these themes in a subsequent paper.

\begin{ack}
The authors are very grateful to Antonio Lei for his interest in the paper and for sharing the preprint \cite{LeiSu}. We would also like to thank Antonio Lei for suggesting a more direct proof of Corollary \ref{main2cor}. We would also like to thank the anonymous referee for providing various helpful comments that have improved the exposition of the paper.
The research of this article took place when S. Ahmed was a postdoctoral fellow at Central China Normal University, and he would like to acknowledge the hospitality
and conducive working conditions provided by the said institute. M. F. Lim gratefully acknowledges support by the
National Natural Science Foundation of China under Grant No. 11550110172 and Grant No. 11771164.
 \end{ack}

\section{Signed Selmer groups} \label{Selmer}

In this section, we will prove Theorem \ref{main1}. As the formula is well-documented when $E$ has good ordinary reduction at every primes of $F$ above $p$ (see \cite[Theorem 3.3]{CS00} or \cite[Theorem 4.1]{G99}), we may and will assume that our elliptic curve $E$ has some primes of supersingular reduction above $p$. In this situation, we shall consider a slightly more general setting following \cite{KO}. As always, $p$ will denote a fixed odd prime. Let $F'$ be a number field and $E$ an elliptic curve defined over $F'$. Fix a finite extension $F$ of $F'$. Let $S$ be a finite set of primes of $F'$ which contains the primes above $p$, the bad reduction primes of $E$, the ramified primes of $F/F'$ and the infinite primes. Denote by $F_S$ the maximal algebraic extension of $F$ which is unramified outside $S$.
For every (possibly infinite) extension $L$ of $F$ contained in $F_S$, we set $G_S(L) = \Gal(F_S/L)$. We shall write $S_p$ (resp., $S_p'$) for the set of primes of $S$ lying above $p$ (resp., not lying above $p$). Denote by $S_p^{ord}$ (resp., $S_p^{ss}$) for the set of good ordinary reduction (resp., good supersingular reduction) primes of $E$ above $p$. We also make the following assumptions:

(S1) The elliptic curve $E$ has good reduction at all primes in $S_p$ and $S_p^{ss}\neq \emptyset$.

(S2) For each $v\in S_p^{ss}$, one has $F'_v=\Qp$ and $a_v = 1 + p - |\tilde{E}_v(\mathbb{F}_p)| = 0$, where $\tilde{E}_v$ is the reduction of $E$ at $v$.

(S3) For each $v\in S_p^{ss}$, $v$ is unramified in $F/F'$.

(S4) For each $w\in S_p^{ss}(F)$, $[F_w:\Qp] \neq 0~(\mathrm{mod}~4)$. Here $S_p^{ss}(F)$ is the set of primes of $F$ above $S_p^{ss}$.

Denote by $F^{\cyc}$ the cyclotomic $\Zp$-extension of $F$ and $F_n$ the intermediate subfield of $F^{\cyc}$ with $|F_n:F|=p^n$. Note that it follows from (S2) and (S3) that every prime $w\in S_p^{ss}(F)$ is totally ramified in $F^{\cyc}/F$. In particular, for each such prime $w$, there is a unique prime of $F_n$ lying above the said prime which we, by abuse of notation, still denote by $w$. Following \cite{KimParity, KimPM, KO, Kob}, we define
the following groups
\[E^+(F_{n,w}) = \{ P\in E(F_{n,w})~:~\mathrm{tr}_{n/m+1}(P)\in E(F_{m,w}), 2\mid m, -1\leq m \leq n-1\}, \]
\[E^-(F_{n,w}) = \{ P\in E(F_{n,w})~:~\mathrm{tr}_{n/m+1}(P)\in E(F_{m,w}), 2\nmid m, -1\leq m \leq n-1\}, \]
where $\mathrm{tr}_{n/m+1}: E(F_{n,w}) \lra E(F_{m+1,w})$ denotes the trace map.

From now on, let $I=\{1,...,r\}$, where $r=|S_p^{ss}(F)|$. We shall index the primes in $S_p^{ss}(F)$ by $w_1,..., w_r$. For each $\overrightarrow{s}=(s_1,...,s_r)\in\{\pm\}^I$, we write \[\mathcal{H}^{\overrightarrow{s}}_n = \bigoplus_{i=1}^r\frac{H^1(F_{n,w_i},E(p))}{E^{s_i}(F_{n,w_i})\ot\Qp/\Zp}.\]
The signed Selmer group is then defined to be
\[\Sel^{\overrightarrow{s}}(E/F_n) = \ker \left(H^1(G_S(F_n),E(p))\lra \mathcal{H}^{\overrightarrow{s}}_n\times\bigoplus_{w\in S_p^{ord}(F_n)}\frac{H^1(F_{n,w},E(p))}{E(F_{n,w})\ot\Qp/\Zp}\times\bigoplus_{w\in S_p'(F_n)}H^1(F_{n,w},E(p)) \right),\]
where $S_p^{ord}(F_n)$ (resp., $S'_p(F_n)$) denotes the set of primes of $F_n$ above $S_p^{ord}$ (resp., $S'_p$).
We also recall that the classical $p$-primary Selmer group for $E$ over $F_n$ is defined by
\[\Sel(E/F_n) = \ker \left(H^1(G_S(F_n),E(p))\lra \bigoplus_{w\in S_p(F_n)}\frac{H^1(F_{n,w},E(p))}{E(F_{n,w})\ot\Qp/\Zp}\times\bigoplus_{w\in S_p'(F_n)}H^1(F_{n,w},E(p)) \right).\]
 The two Selmer groups fit into the following commutative diagram
\[   \xymatrixrowsep{0.25in}
\xymatrixcolsep{0.15in}\entrymodifiers={!! <0pt, .8ex>+} \SelectTips{eu}{}\xymatrix{
    0 \ar[r]^{} & \Sel^{\overrightarrow{s}}(E/F_n) \ar[d]_{\al} \ar[r] &  H^1(G_S(F_n), E(p))
    \ar@{=}[d] \ar[r]^(.27){\psi^{\overrightarrow{s}}} & \displaystyle\mathcal{H}^{\overrightarrow{s}}_n\times\bigoplus_{w\in S_p^{ord}(F_n)}\frac{H^1(F_{n,w},E(p))}{E(F_{n,w})\ot\Qp/\Zp}\times\bigoplus_{w\in S_p'(F_n)}H^1(F_{n,w},E(p))\ar[d]\\
    0 \ar[r]^{} & \Sel(E/F_{n}) \ar[r]^{} & H^1(G_S(F_n), E(p)) \ar[r]^(.3){\phi} & \displaystyle\bigoplus_{w\in S_p(F_n)}\frac{H^1(F_{n,w},E(p))}{E(F_{n,w})\ot\Qp/\Zp}\times\bigoplus_{w\in S_p'(F_n)}H^1(F_{n,w},E(p)) } \]
with exact rows. Denote by $\psi^{\overrightarrow{s}}_{ss}$ the map from $\Sel(E/F_n)$ to $\mathcal{H}^{\overrightarrow{s}}_n$ that is induced by $\psi^{\overrightarrow{s}}$. It is now straightforward to verify the following.

\bl \label{PM Selmer} We have the following identification
\[\Sel^{\overrightarrow{s}}(E/F_n) = \ker \left(\Sel(E/F_n)\stackrel{\psi^{\overrightarrow{s}}_{ss}}{\lra} \mathcal{H}^{\overrightarrow{s}}_n\right).\]
\el

Write $\Sel^{\overrightarrow{s}}(E/F^{\cyc})=\ilim_n \Sel^{\overrightarrow{s}}(E/F_n)$ and $\mathcal{H}^{\overrightarrow{s}}_{\infty}= \ilim_n \mathcal{H}^{\overrightarrow{s}}_n$. It is not difficult to verify that $\Sel^{\overrightarrow{s}}(E/F^{\cyc})$ is cofinitely generated over $\Zp\ps{\Ga}$. In fact, one expects the following conjecture which is a natural extension of Mazur \cite{Maz} and Kobayashi \cite{Kob}.

\medskip \noindent \textbf{Conjecture.} $\Sel^{\overrightarrow{s}}(E/F^{\cyc})$ is a cotorsion $\Zp\ps{\Ga}$-module, where $\Ga=\Gal(F^{\cyc}/F)$.

\medskip
When $S_p^{ss}$ is empty, the above conjecture is precisely Mazur's conjecture \cite{Maz} which is known in the case when $E$ is defined over $\Q$ and $F$ an abelian extension of $\Q$ (see \cite{K}). When $E$ is an elliptic curve over $\Q$ with good supersingular singular reduction at $p$, this conjecture was established by Kobayashi
(cf. \cite{Kob}; also see \cite{BL} for some recent progress on this conjecture). Here we shall prove the following. Theorem \ref{main1} will follow from this by taking $F=F'$.

\bt \label{main2}
Assume that $(S1)-(S4)$ are valid. Suppose that $\Sel(E/F)$ is finite. Then $\Sel^{\overrightarrow{s}}(E/F^{\cyc})$ is a cotorsion $\Zp\ps{\Ga}$-module and its
 $\Ga$-Euler characteristics is given by \[ \frac{|\Sha(E/F)(p)|}{|E(F)(p)|^2}\times\prod_w c_w^{(p)}\times\prod_{w\in S_p^{ord}(F)}(d_w^{(p)})^2.\]
\et

The remainder of this section will be devoted the proof of Theorem \ref{main2}. As a start, we record the following two preparatory lemmas which are required for our calculation.

\bl \label{E(p)=0} Assume that $(S1)-(S3)$ are valid. Then $E(F)(p)=0$ and $E(F^{\cyc})(p)=0$.
\el

\bpf
For $w\in S_p^{ss}(F)$, a similar argument to that in \cite[Proposition 8.7]{Kob} tells us that $E(F_w)(p)=0$. Since we are assuming that $S_p^{ss}\neq \emptyset$, this in turn implies that $E(F)(p)=0$. But as $F^{\cyc}/F$ is a pro-$p$ extension, it follows from \cite[Corollary 1.6.13]{NSW} that $E(F^{\cyc})(p)=0$.
\epf

\bl \label{H2=0} Assume that $(S1)-(S3)$ are valid. Suppose that $\Sel(E/F)$ is finite. Then we have that $H^2\big(G_S(F^{\cyc}),E(p)\big)=0$,  $H^1\big(\Ga, H^1(G_S(F^{\cyc}),E(p))\big)=0$ and
\[ H^1(G_S(F),E(p))\cong H^1(G_S(F^{\cyc}),E(p))^{\Ga}.\]
\el

\bpf
Since $\Ga$ has $p$-cohomological dimension one, the spectral sequence
\[ H^i\big(\Ga,H^j(G_S(F^{\cyc}),E(p))\big)\Longrightarrow H^{i+j}(G_S(F),E(p))\]
yields short exact sequences
\[ 0\lra H^1\big(\Ga, E(F^{\cyc})(p))\big)\lra H^1(G_S(F),E(p))\lra H^1(G_S(F^{\cyc}),E(p))^{\Ga}\lra 0 \] and
\[ 0\lra H^1\big(\Ga, H^1(G_S(F^{\cyc}),E(p))\big)\lra H^2(G_S(F),E(p))\lra H^2(G_S(F^{\cyc}),E(p))^{\Ga}\lra 0. \]
The final isomorphism of the lemma follows from the first short exact sequence and Lemma \ref{E(p)=0}.
On the other hand, as $\Sel(E/F)$ is finite, it follows from \cite[Proposition 1.9]{CS00} that $H^2(G_S(F),E(p))=0$. Putting this into the second short exact sequence, we obtain $H^1\big(\Ga, H^1(G_S(F^{\cyc}),E(p))\big)=0$ and $H^2(G_S(F^{\cyc}),E(p))^{\Ga}=0$, where the latter in turn implies that $H^2(G_S(F^{\cyc}),E(p))=0$.
This proves the lemma.
\epf

The next lemma is concerned with analysing the local map
\[g_w : \frac{H^1(F_w, E(p))}{E(F_v)\ot \Qp/\Zp} \lra \left(\frac{H^1(F^{\cyc}_w, E(p))}{E^{\pm}(F^{\cyc}_w)\ot \Qp/\Zp}\right)^{\Ga} \]
for $w\in S_p^{ss}(F)$.

\bl \label{local calculation} Suppose that $(S1)-(S4)$ are valid.  Then for every $w\in S_p^{ss}(F)$, the map $g_w$ is an isomorphism.
\el

\bpf
We essentially follow the idea in the proof of \cite[Proposition 4.28]{KimParity}.
Consider the following diagram
\[   \entrymodifiers={!! <0pt, .8ex>+} \SelectTips{eu}{}\xymatrix{
    0 \ar[r]^{} & E(F_v)\ot \Qp/\Zp \ar[d]_{a_w} \ar[r] &  H^1(F_w, E(p))
    \ar[d]_{b_w} \ar[r] & \displaystyle\frac{H^1(F_w, E(p))}{E(F_v)\ot \Qp/\Zp} \ar[d]_{g_w} \ar[r]^{} & 0\\
    0 \ar[r]^{} & \big(E^{\pm}(F^{\cyc}_w)\ot \Qp/\Zp\big)^{\Ga} \ar[r]^{} & H^1(F^{\cyc}_w, E(p))^{\Ga} \ar[r] & \displaystyle\left(\frac{H^1(F^{\cyc}_w, E(p))}{E^{\pm}(F^{\cyc}_w)\ot \Qp/\Zp}\right)^{\Ga} &} \]
    with exact rows. As seen from the proof of Lemma \ref{E(p)=0}, we have $E(F_w)(p)=0$ which in turn implies that $E(F^{\cyc}_w)(p)=0$. Hence we have that $b_w$ is an isomorphism. Consequently, $a_w$ is injective. By (S4) and \cite[Corollary 3.25]{KO}, we have that $\big(E^{\pm}(F^{\cyc}_w)\ot \Qp/\Zp\big)^{\Ga}$ is a cofree $\Zp$-module with $\Zp$-corank $[F_w:\Qp]$. But by Mattuck's theorem \cite{Mat}, $E(F_v)\ot \Qp/\Zp$ is also a cofree $\Zp$-module with $\Zp$-corank $[F_w:\Qp]$. Hence $a_w$ has to be an isomorphism which in turn implies that $g_w$ is injective.

    But now, using the fact that $E(F_w)(p)=0$, it then follows from local Tate duality that $H^2(F_w, E[p]) =0$ which in turn implies that $H^1(F_w, E(p))$ is $p$-divisible. Combining this latter observation with a standard local Euler characteristic calculation (cf. \cite[\S 3, Proposition 1]{G89}), we have that $H^1(F_w, E(p))$ is a cofree $\Zp$-module with $\Zp$-corank $[F_w:\Qp]$ . On the other hand, it follows from \cite[Proposition 3.32]{KO} that $\displaystyle\left(\frac{H^1(F^{\cyc}_w, E(p))}{E^{\pm}(F^{\cyc}_w)\ot \Qp/\Zp}\right)^{\Ga}$ is a cofree $\Zp$-module with $\Zp$-corank $[F_w:\Qp]$. Thus, we have shown that $g_w$ is an injection between two $p$-divisible groups of the same $\Zp$-corank and hence it must be an isomorphism. This proves the lemma.
\epf

Now consider the following diagram
\[   \entrymodifiers={!! <0pt, .8ex>+} \SelectTips{eu}{}\xymatrix{
    0 \ar[r]^{} & \Sel(E/F) \ar[d]_{a} \ar[r] &  H^1(G_S(F), E(p))
    \ar[d]_{h} \ar[r]^(.35){\rho} & \displaystyle\displaystyle\bigoplus_{w|p}\frac{H^1(F_w,E(p))}{E(F_{v})\ot\Qp/\Zp}\times\bigoplus_{w\in
    S_p'(F)}H^1(F_w,E(p))\ar[d]_{g=\oplus_w g_w}\\
    0 \ar[r]^{} & \Sel^{\overrightarrow{s}}(E/F^{\cyc})^{\Ga} \ar[r]^{} & H^1(G_S(F^{\cyc}), E(p))^{\Ga} \ar[r]^(.4){\phi_{\infty}} & \left(\displaystyle\mathcal{H}^{\overrightarrow{s}}_{\infty}\times
    \mathcal{H}^{ord}_{\infty}\times\bigoplus_{w\in S_p'(F^{\cyc})}H^1(F^{\cyc}_w,E(p))\right)^{\Ga} } \]
 with exact rows, where $\mathcal{H}^{ord}_{\infty}=\displaystyle\ilim_n\bigoplus_{w\in S_p^{ord}(F_n)}\frac{H^1(F_{n,w},E(p))}{E(F_{n,w})\ot\Qp/\Zp}$. We shall make use of the notation in the above diagram without further mention.

\bl \label{H1=0} Assume that $(S1)-(S4)$ are valid. Suppose that $\Sel(E/F)$ is finite. Then $\rho$ is surjective and $H^1(\Ga, \Sel^{\overrightarrow{s}}(E/F^{\cyc}))=0$.
\el

\bpf
Since $\Sel(E/F)$ is finite, it follows from \cite[Proposition 1.9]{CS00} that $\coker \rho$ is finite of order $|E(F)(p)|$. By Lemma \ref{E(p)=0}, this in turn implies that $\rho$ is surjective which proves the first assertion of the lemma.

Combining \cite[Lemma 3.4 and Proposition 3.5]{CS00} with Lemma \ref{local calculation}, we have that $g$ is surjective. Therefore, $\phi_{\infty}$ is also surjective. Now consider the following exact sequence \[
       0 \lra \Sel^{\overrightarrow{s}}(E/F^{\cyc}) \lra H^1(G_S(F^{\cyc}), E(p)) \stackrel{\phi}{\lra}  B,  \]
 where $B= \displaystyle\mathcal{H}^{\overrightarrow{s}}_{\infty}\times
       \mathcal{H}^{ord}_{\infty}\times\bigoplus_{w\nmid p}H^1(F^{\cyc}_w,E(p))$. Write $A= \mathrm{im}(\phi)$ and $C= \coker (\phi)$.
Taking $\Ga$-invariant of the short exact sequence
\[ 0\lra \Sel^{\overrightarrow{s}}(E/F^{\cyc}) \lra H^1(G_S(F^{\cyc}), E(p)) \lra A\lra 0,\]
and taking Lemma \ref{H2=0} into account,
we obtain an exact sequence
\[ 0\lra \Sel^{\overrightarrow{s}}(E/F^{\cyc})^{\Ga} \lra H^1(G_S(F^{\cyc}), E(p))^{\Ga} \stackrel{\tau}\lra A^{\Ga}\lra H^1(\Ga, \Sel^{\overrightarrow{s}}(E/F^{\cyc}))\lra 0\]
with $H^1(\Ga, A)=0$. Taking the latter into consideration, it follows from
the $\Ga$-invariant of the short exact sequence
\[ 0\lra A \lra B \lra C\lra 0 \] that we obtain
a short exact sequence
\[ 0\lra A^{\Ga} \lra B^{\Ga} \lra C^{\Ga}\lra 0.\]
Since $\phi_{\infty}$ is surjective and it is given by the composition  $H^1(G_S(F^{\cyc}),E(p))\lra A^{\Ga} \lra B^{\Ga}$, we have that the injection $A^{\Ga} \lra B^{\Ga}$ is also surjective and hence an isomorphism. Under this identification, we have $\tau =\phi_{\infty}$, whose surjectivity in turn implies that $H^1(\Ga, \Sel^{\overrightarrow{s}}(E/F^{\cyc}))=0$. The proof of the lemma is now completed.\epf

We record the following by-product of our argument which is not required for the final proof. It may also be quite possible that one can derive the conclusion of this said result via the methods of \cite[Proposition 3.10]{KimPM}. However, we decide to include the following alternative proof which might be of interest in its own right. We should however mention that our proof here relies on the finiteness assumption of $\Sel(E/F)$.

\bp \label{surjective} Assume that $(S1)-(S4)$ are valid. Suppose that $\Sel(E/F)$ is finite. Then we have the following short exact sequence \[
       0 \lra \Sel^{\overrightarrow{s}}(E/F^{\cyc}) \lra H^1(G_S(F^{\cyc}), E(p)) \stackrel{\phi}{\lra}  \displaystyle\mathcal{H}^{\overrightarrow{s}}_{\infty}\times
       \mathcal{H}^{ord}_{\infty}\times\bigoplus_{w\in S_p'(F^{\cyc})}H^1(F^{\cyc}_w,E(p))\lra 0.  \]
\ep

\bpf
We retain the notation of Lemma \ref{H1=0}. From the proof of the said lemma, we have obtained
a short exact sequence
\[ 0\lra A^{\Ga} \lra B^{\Ga} \lra C^{\Ga}\lra 0\]
and shown that  $A^{\Ga} \cong B^{\Ga}$. Thus, we have $C^{\Ga}=0$ which in turn implies that $C=0$. But recall that $C= \coker \phi$ and so this proves the proposition.
\epf

We can finally prove Theorem \ref{main2}.

\bpf[Proof of Theorem \ref{main2}] To prove the first assertion of the theorem, it suffices to show that $\Sel^{\overrightarrow{s}}(E/F^{\cyc})^{\Ga}$ is finite. By Lemma \ref{H2=0}, $h$ is an isomorphism. Therefore, by the snake lemma, we are reduced to showing that $\ker g$ is finite. In fact, for $w\in S_p^{ord}(F)$, $\ker g_w$ is finite with order $(d_w^{(p)})^2$ (cf. \cite[Proposition 3.5]{CS00} or \cite[Lemma 4.4]{G99}).  If $w\in S_p^{ss}(F)$, $g_w$ is an isomorphism by Lemma \ref{local calculation}. Finally, for $w\nmid p$, $\ker g_w$ is finite with order $c_w^{(p)}$ (cf. \cite[Lemma 3.4]{CS00} or \cite[Lemma 4.4]{G99}). Hence $\ker g$ is finite as required.

It remains to compute the $\Ga$-Euler characteristics of $\Sel^{\overrightarrow{s}}(E/F^{\cyc})$. By Lemma \ref{H2=0}, $\rho$ is surjective. Taking the final isomorphism in the assertion of Lemma \ref{H2=0} into account, it then follows from the above diagram that
\[ \left|\Sel^{\overrightarrow{s}}(E/F^{\cyc})^{\Ga}\right| = |\Sel(E/F)|~|\ker g|.\] By Lemma \ref{H1=0}, the left hand side is precisely the $\Ga$-Euler characteristics of
$\Sel^{\overrightarrow{s}}(E/F^{\cyc})$. Since $\Sel(E/F)$ is finite, we have $|\Sel(E/F)| = |\Sha(E/F)(p)|$. Also, as seen above, we have that $|\ker g|$ is given by  $\prod_w c_w^{(p)}\times\prod_{w\in S_p^{ord}(F)}(d_w^{(p)})^2$. Combining these calculations, we obtain the required formula noting that $|E(F)(p)|=1$ by Lemma \ref{E(p)=0}.
\epf

We record an interesting corollary of (the proof of) our Theorem \ref{main2}.

\bc \label{main2cor}
Assume that $(S1)-(S4)$ are valid. Suppose that there exists $\overrightarrow{t}\in \{\pm\}^{I}$ such that $\Sel^{\overrightarrow{t}}(E/F^{\cyc}) =0$. Then $\Sel^{\overrightarrow{s}}(E/F^{\cyc})=0$ for every $\overrightarrow{s}\in \{\pm\}^{I}$.
\ec

\bpf
Suppose that $\Sel^{\overrightarrow{t}}(E/F^{\cyc}) =0$ for some $\overrightarrow{t}\in \{\pm\}^{I}$. Then from the diagram before Lemma \ref{surjective}, we have that $\Sel(E/F)=0$. In particular, $\Sel(E/F)$ is finite. Therefore, we apply the argument of the proof of Theorem \ref{main2} to obtain the equality
\[ \left|\Sel^{\overrightarrow{t}}(E/F^{\cyc})^{\Ga}\right| = |\Sel(E/F)|~|\ker g| = |\ker g|.\]
Since $\Sel^{\overrightarrow{t}}(E/F^{\cyc}) =0$, it follows that $\ker g =0$. From the proof of Theorem \ref{main2}, we also see that $\ker g$  have the same common value for every $\overrightarrow{s}\in \{\pm\}^{I}$ and hence is trivial. Consequently, we have
\[\left|\Sel^{\overrightarrow{s}}(E/F^{\cyc})^{\Ga}\right| =0\]
which in turn implies that $\Sel^{\overrightarrow{s}}(E/F^{\cyc})^{\Ga}=0$. The latter is of course equivalent to saying that $\Sel^{\overrightarrow{s}}(E/F^{\cyc})=0$ as required.
\epf

\section{Concluding remarks} \label{concluding remarks}

We make some remarks about Theorem \ref{main2}. In this said theorem, we have made an assumption that for each $w\in S_p^{ss}(F)$, $[F_w:\Qp] \neq 0~(\mathrm{mod}~4)$ (this is our assumption (S4)). We should mention that if all the signs appearing in the signed Selmer group are $-$, one does not require this assumption (S4). However, if at least one of the signs is a $+$, , we will not be able to prove that the local map $g_w$ is injective without the said assumption. In fact, tracing the proof of Lemma \ref{local calculation}, it would seem that $g_w$ has kernel which is a cofree $\Zp$-module with corank $2$ (when $[F_w:\Qp] = 0~(\mathrm{mod}~4)$). This seems reminiscent of the so-called ``exceptional zeroes'' phenomenon in the case of a split multiplicative prime (for instances, see \cite{G94, MazTT}). We do not have a good explanation on this at this point of writing but we hope to come back to this issue in a future work.

\end{document}